# An Extended Admittance Modeling Method with Synchronization Node for Stability Assessment of Converters-Interlinked System

Haoxiang Zong, *Member, IEEE*, Chen Zhang, *Member, IEEE*, and Marta Molinas, *Fellow, IEEE*

***Abstract*—Diverse synchronization dynamics within the grid-following (GFL)/grid-forming (GFM) converters-interlinked system are prone to induce oscillatory instabilities. To quantify their stability influences, frequency-domain modal analysis (FMA) method based on the impedance network can serve as a good reference. However, since the adopted impedance network only retains electrical nodes, oscillation information provided by the FMA method is mainly concerned with circuits (e.g., participation of nodes), which is not convenient for an intuitive probe of sync loops' participations. To address this issue, this paper proposes an extended admittance modeling method for FMA, the basis of which is the explicit characterization of GFL/GFM sync loops. First, a four-port extended impedance model (EIM) of converter with one virtual sync node is proposed. Its resulting extended impedance network (EIN) is formed for the converters-interlinked system. Then, the FMA method can be directly applied to those virtual sync nodes/branches, so as to realize an intuitive evaluation of sync dynamics' effects on oscillations. The effectiveness of the proposed method is validated by the frequency scanning and time domain simulations in a typical point-to-point HVDC system.**

## I. Introduction

LATELY, a steadily increasing share of converter-interfaced renewable energy sources (RPGs) along with a decreased ratio of synchronous generators has led to profound operation issues of the power system [1], e.g., low inertia and oscillatory instability, etc. One tentative solution of great promise is to duly replace a certain proportion of grid-following (GFL) converters (usually as grid-connection interfaces of RPGs) with grid-forming (GFM) ones [2]. Due this, the hybrid or heterogeneous power system composed of GFL/GFM converters is becoming more common, typical as the wind farm with a battery energy storage system (BESS) [3] in which wind turbines operate at GFL mode while BESS is GFM controlled for the enhancement of the grid adaptability.

A salient feature of hybrid GFL/GFM converter systems is the coexistence of diverse synchronization (briefly denoted as "sync") behaviors, including phase-locked loop (PLL), droop control, virtual synchronous generator (VSG) control, matching control [4], etc. [5] shows that there is a great complementarity between GFL and GFM converters in terms of grid integration characteristics (e.g., autonomous inertia response of GFM versus fast power response of GFL), and thus a combination of both is expected to be beneficial to system operations. However, control effects of GFL/GFM converters are nonideal in practice, where control interactions among interlinked converters are readily provoked, further leading to oscillatory instabilities. Recent studies have shown that phase-locked loops (PLL) and their interactions play a crucial role in the grid-synchronization stability of GFL systems [6], [7]. Similar synchronization stability issues also exist in GFM converters [8] according to the well-known swing dynamics of SGs. Therefore, it can be anticipated that sync dynamics induced oscillation issues will be more evident in hybrid GFL/GFM systems [9], which are of great value for analysis. In this regard, an appropriate modeling and analysis approach is required, particularly be capable of explicitly characterizing diverse sync dynamics while at the same time preserving good scalability for the system-level analysis. This is also the main objective of this paper.

Regarding the oscillatory stability analysis of multi converters -interlinked systems, the frequency-domain modal analysis (FMA) method based on the system's impedance network characterization can be applied and is drawing great attention in recent years [10]-[15]. The application of FMA mainly consists of three steps: 1) characterize electric port dynamics of each apparatus as circuit impedances, including *ac* [16], *dc* [17], as well as three-port *ac*/*dc* [18] impedances; 2) formulate the systemwide impedance network by assembling impedances of all components according to the network topology [19]-[24]; 3) calculate sensitivity indices like participation factor (PF) of nodes, branches, components, etc. With these quantitative indices, instability dominant components and oscillation propagation behaviors can be evaluated. It should be mentioned that FMA method can be applied to long-transmission line systems [10], ac-interconnected systems [11]-[13] and hybrid AC/DC system [14], [15], by adopting suitable *ac* [19]-[21], or *ac*/*dc* [22]-[24] impedance network models.

This work was supported by the National Natural Science Foundation of China under Grant 52207215. *(Corresponding author: Chen Zhang).*

Haoxiang Zong and Chen Zhang are with the Key Laboratory of Control of Power Transmission and Conversion of Ministry of Education, Shanghai Jiao Tong University (SJTU), 200240, Shanghai, China, also with the Department of Electrical Engineering, School of Electronic Information and Electrical Engineering, SJTU, 200240, Shanghai. (e-mails: {haoxiangzong, nealbc} @sjtu.edu.cn ).

Marta Molinas is with Department of Engineering Cybernetics, Norwegian University of Science and Technology, Trondheim, 7491, Norway. (e-mail: marta.molinas@ntnu.no).

However, conventional impedance modeling techniques as mentioned above only retain electrical ports [21], making the application of FMA method being limited to the acquisition of oscillation properties concerning circuits (e.g., participation of nodes, branches, [11] etc.). This implies that properties of oscillatory instability arising from GFL/GFM sync loops cannot be intuitively probed. A trade-off approach is to conduct the parameter sensitivity analysis [13], but this requires a complex symbolic computation manipulated by multi-layer partial derivatives (usually from nodes to components and finally into parameters). Furthermore, entries of system closed-loop matrix that contain sync control parameters are hard to be predefined, meaning each matrix element should be went through and implemented with partial derivative calculus [15]. Therefore, the participation identification and analysis of sync loops are still challenging. To this end and oriented for better reflection of sync loops' interactions within hybrid GFL/GFM converters-interlinked systems, this paper proposes a sync-node extended modeling method, which has following innovative aspects and is expected to contribute to the general aspect of the development of the impedance-based FMA method:

1) An extended and modularized impedance representation of a generic ac/dc converter able to explicitly characterize the sync loop is proposed, referred to as the four-port Extended Impedance Model (EIM).
2) A method to formulate the Extended Impedance Network (EIN) using the four-port EIM is further established, with which a sync-node extended modal analysis for better displaying the sync contribution is demonstrated.
3) Interactions between typical GFM and GFL sync loops on the system oscillatory stability are revealed.

The rest of paper is arranged as: Section II briefly reviews the current FMA method. Section III introduces the sync-node extended FMA method. Section IV demonstrates the feasibility and capacity of the proposed method by case studies of a point-to-point HVDC system. Section V concludes the paper.

## II. Review on the Conventional FMA Method

A brief review of the current FMA method [10]-[15] is given to better introduce proposition of this work. As aforementioned, the typical FMA method is achieved by performing frequency domain modeling along with sensitivity analysis. Regarding the modeling part, two-port ac admittance and three-port ac/dc admittance as depicted by (1) are often adopted.

$$\text{two port}: \Delta \boldsymbol{i}_{gdq}(s) = \boldsymbol{Y}_{dq}(s)\Delta \boldsymbol{u}_{gdq}(s)$$

$$\text{three port}: \begin{bmatrix} \Delta \boldsymbol{i}_{gdq}(s) \\ \Delta i_{dc}(s) \end{bmatrix} = \begin{bmatrix} \boldsymbol{Y}_{dq}(s) & \boldsymbol{a}_{dq}(s) \\ \boldsymbol{b}_{dq}(s) & Y_{dc}(s) \end{bmatrix} \begin{bmatrix} \Delta \boldsymbol{u}_{gdq}(s) \\ \Delta u_{dc}(s) \end{bmatrix} \quad (1)$$

Both of them are acquired generally by first linearizing the converter model with controls, then evaluating the input-output relation of currents $\Delta \boldsymbol{i}_{gdq}$, $\Delta i_{dc}$ and voltages $\Delta \boldsymbol{u}_{gdq}$, $\Delta u_{dc}$ at converter's ac- and dc-side. The two-port ac impedance is usually used for pure ac systems [19], while the three-port ac/dc admittance is feasible for hybrid ac/dc systems, e.g., a MTDC system [23]. Further, given impedances of each component are assembled into a circuit network as in Fig. 1(a), by partitioning the overall system into the converter- and network-subsystem at the point-of-connections (PoCs) and applying impedance-based modeling methods [24], loop-gain $\boldsymbol{L}(s)$ characterization of the system can be obtained as (2), of which the closed-loop stability can be judged by zeros of $\det[\boldsymbol{I}+\boldsymbol{L}(s)]=0$.

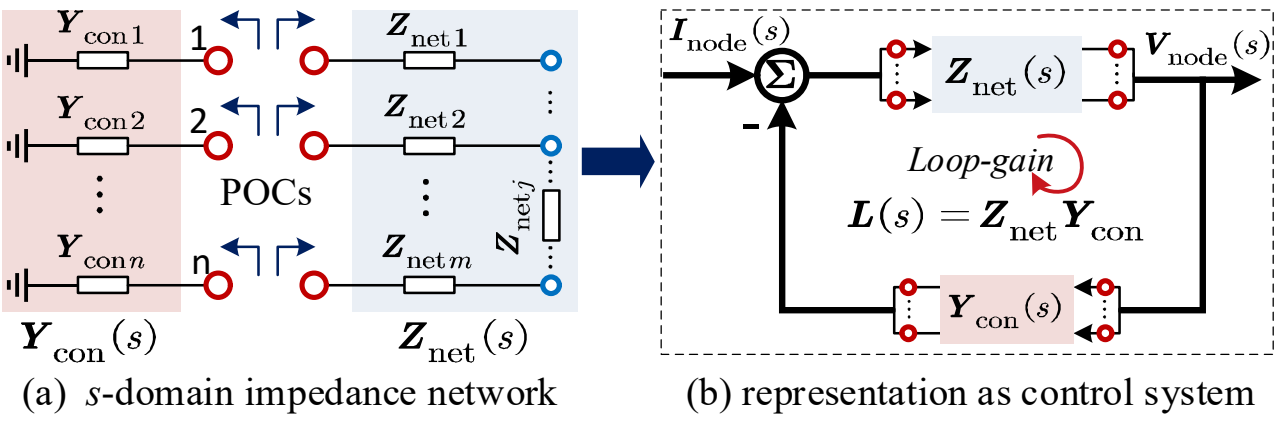

Fig. 1. $\boldsymbol{L}(s)$-based representations of a converters-dominated system. "con" denotes "converter", and "net" denotes passive electric "networks".

$$\boldsymbol{V}_{node}(s) = [\boldsymbol{I} + \underbrace{\boldsymbol{Z}_{net}(s)\boldsymbol{Y}_{con}(s)}_{\boldsymbol{L}(s)}]^{-1}\boldsymbol{Z}_{net}(s)\cdot \boldsymbol{I}_{node}(s) \quad (2)$$

To gain more information on oscillatory stability properties, sensitivity analysis using eigen-decomposition information (see (3)) of $\boldsymbol{L}$ can be further performed, which is the typical FMA [10]. In detail, suppose one mode $\lambda_m$ is interested for analysis, by evaluating the eigenvalues of $\boldsymbol{I}+\boldsymbol{L}(\lambda_m)$, an entry that is zero can be located, e.g., the $k$-th one denoted as $1+\Lambda_k$. Then, the frequency-domain PF is defined as the sensitivity of $1+\Lambda_k$ to an entry in $\boldsymbol{I}+\boldsymbol{L}$, which can be expressed as (4).

$$\boldsymbol{L} = \boldsymbol{R}\boldsymbol{\Lambda}\boldsymbol{T} \quad (3)$$

$$\frac{\partial 1+\Lambda_k}{\partial L_{ij}} = \frac{\partial \Lambda_k}{\partial L_{ij}} = \boldsymbol{t}_k \frac{\partial \boldsymbol{L}}{\partial L_{ij}} \boldsymbol{r}_k = \left[\boldsymbol{r}_k \boldsymbol{t}_k\right]_{ji} = \left[\boldsymbol{PF}\right]_{ji} \quad (4)$$

where $\Lambda=\text{diag}(\Lambda_1,\Lambda_2,\ldots,\Lambda_n)$ is the eigenvalue matrix of $\boldsymbol{L}$; $\boldsymbol{R}=[\boldsymbol{r}_1,\boldsymbol{r}_2,\ldots,\boldsymbol{r}_n]$, $\boldsymbol{T}=[(\boldsymbol{t}_1)^T,(\boldsymbol{t}_2)^T,\ldots,(\boldsymbol{t}_n)^T]^T$ are right and left eigenvectors. If diagonal entries of $\boldsymbol{L}$ are analyzed, the resulting $\boldsymbol{PF}$ can reflect the sensitivity of principal loop (i.e., $\boldsymbol{Z}_{net\,i}\boldsymbol{Y}_{con\,i}$) at each node/PoC to the critical loop gain $\Lambda_k$.

With these sensitivity indices, significance of electric nodes to which either converters or passive components are connected on system instability can be identified. Such information can be helpful in synthesizing oscillation attenuation countermeasures, e.g., active damping placement. However, this conventional FMA method is in the sense of pure circuit analysis, while effects of sync control loops can hardly be identified (as they are implicitly distributed inside elements of $\boldsymbol{L}$, i.e., $\boldsymbol{Y}_{con\,i}(s)$). To address this issue, it will be best to explicitly model various sync loops of converter, and it is shown next that how this is achieved in the impedance modeling framework so that the scalability of impedance approach can be preserved.

## III. Sync Node-Extended FMA Method

This section presents a sync-node extended modeling method for the promotion of the FMA method to sync related studies. The idea is that sync loops are separated out from the converter control and explicitly augmented as an external sync port to the existing three-port impedance model like (1), resulting in the *four-port* EIM of this work; then, by regarding the added sync port as a virtual electric path, the EIN consisting of EIMs can be formulated with typical networked formulation method, e.g., node admittance method; finally, the EIN-based FMA is performed, so that sensitivity analysis with sync ports can be better achieved. These parts are described below.

*A. Formulation of the Four-Port EIM for VSC*

As shown in Fig. 2, the typical three-port impedance of VSC is depicted by Layer-1, which consists of electric ports 1~3 and the sync loop is embedded inside the complete control system. To explicitly characterize sync dynamics, the sync loop is separated from the control system and augmented as a virtual port (i.e., port 4) as shown in Fig. 2. Under this modification, the converter is characterized by a four-port module as shown in Layer2, which is referred as the EIM in this paper. In this way, physical structure of the sync control loop is kept intact, and sync dynamics can be studied explicitly. The establishment of this four port EIM is described below.

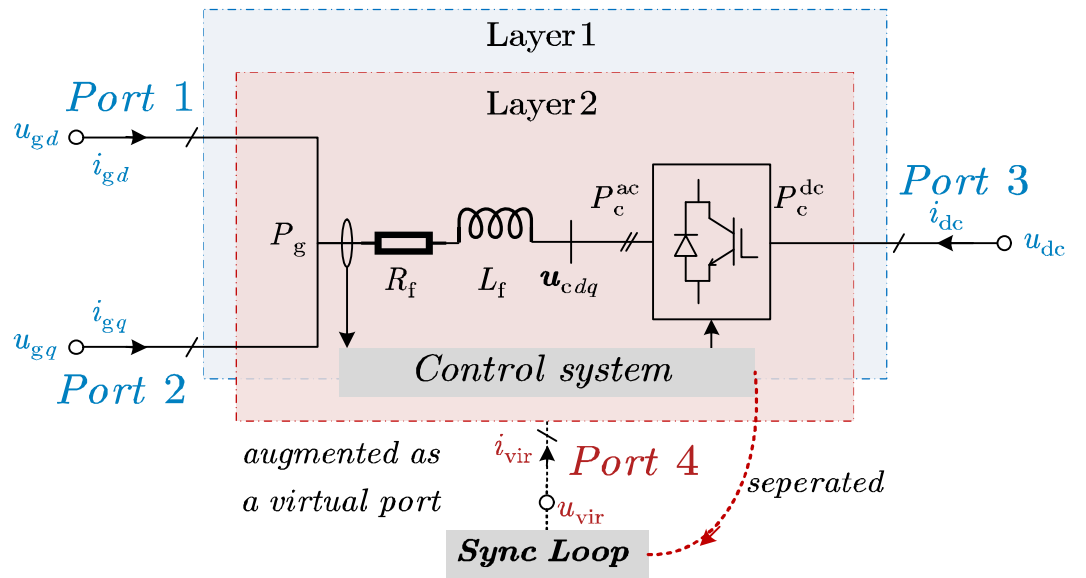

Fig. 2. Schematic diagram of the sync port separation.

*1) Unified Representation of GFL/GFM Sync Loops*

To modularize the four port EIM, a unified representation of separated GFL/GFM sync loops is needed. Suitably, as shown in Fig. 3(a), GFL and GFM sync loops can be represented in a unified feedback form. The sync forward path $G(s)$ consists of a PI regulator and its output is the sync frequency $\omega_{\mathrm{sync}}$. The sync feedback path is constituted of an integrator and the respective feedback transfer function $F(s)$, the output of which could be the active power in VSG, $q$-axis voltage in PLL or dc voltage in matching control [4]. Details of commonly-used sync control method of the GFL and GFM are given in Fig. 3(b).

According to above description, small-signal characteristic of the unified forward path of the sync loop can be written as:

$$\underset{\Delta\omega_{\mathrm{sync}}(s)}{\overset{\Delta u_{\mathrm{vir}}}{\downarrow}} = \underset{Z_{\mathrm{sync}}^{\mathrm{fo}}(s)}{\overset{G(s)}{\downarrow}} \cdot \underset{\Delta P_{\mathrm{sync}}(s)}{\overset{\Delta i_{\mathrm{vir}}}{\downarrow}} \tag{5}$$

where the input $\Delta P_{\mathrm{sync}}$ and output $\Delta\omega_{\mathrm{sync}}$ of the sync loop can be regarded as the *virtual current* $i_{\mathrm{vir}}$ and *virtual voltage* $u_{\mathrm{vir}}$ (as signified in Fig. 2). Under this emulated convention, $G(s)=H_{\mathrm{pll}}(s)$ or $H_{\mathrm{vsg}}(s)$ can be seen as a kind of virtual impedance, thus denoted by symbol $Z_{\mathrm{sync}}^{\mathrm{fo}}$. To keep the unified form, the denotation $\Delta P_{\mathrm{sync}}$ is adopted for both PLL and VSG, where $\Delta P_{\mathrm{sync}}$ in PLL refers to the $q$-axis voltage $\Delta u_{gq}^{c}$ in controller frame and $\Delta P_{\mathrm{sync}}$ in VSG refers to the active power $-\Delta P_{\mathrm{g}}$. The meaning of $\Delta P_{\mathrm{sync}}$ can be altered according to its corresponding sync loop.

Sync dynamics of converters are dominated by the above sync forward path, for instance, inertia constant $J$ and damping ratio $D$ in $H_{\mathrm{vsg}}(s)$ of VSG's sync forward path affect sync behaviors most. In view of this, sync forward path $Z_{\mathrm{sync}}^{\mathrm{fo}}$ is categorized as the external component just as those AC/DC grids connected to converter. As to the sync feedback path shown in Fig. 3(b), it will be integrated into the following four-port EIM, because it is coupled with other control dynamics of the converter.

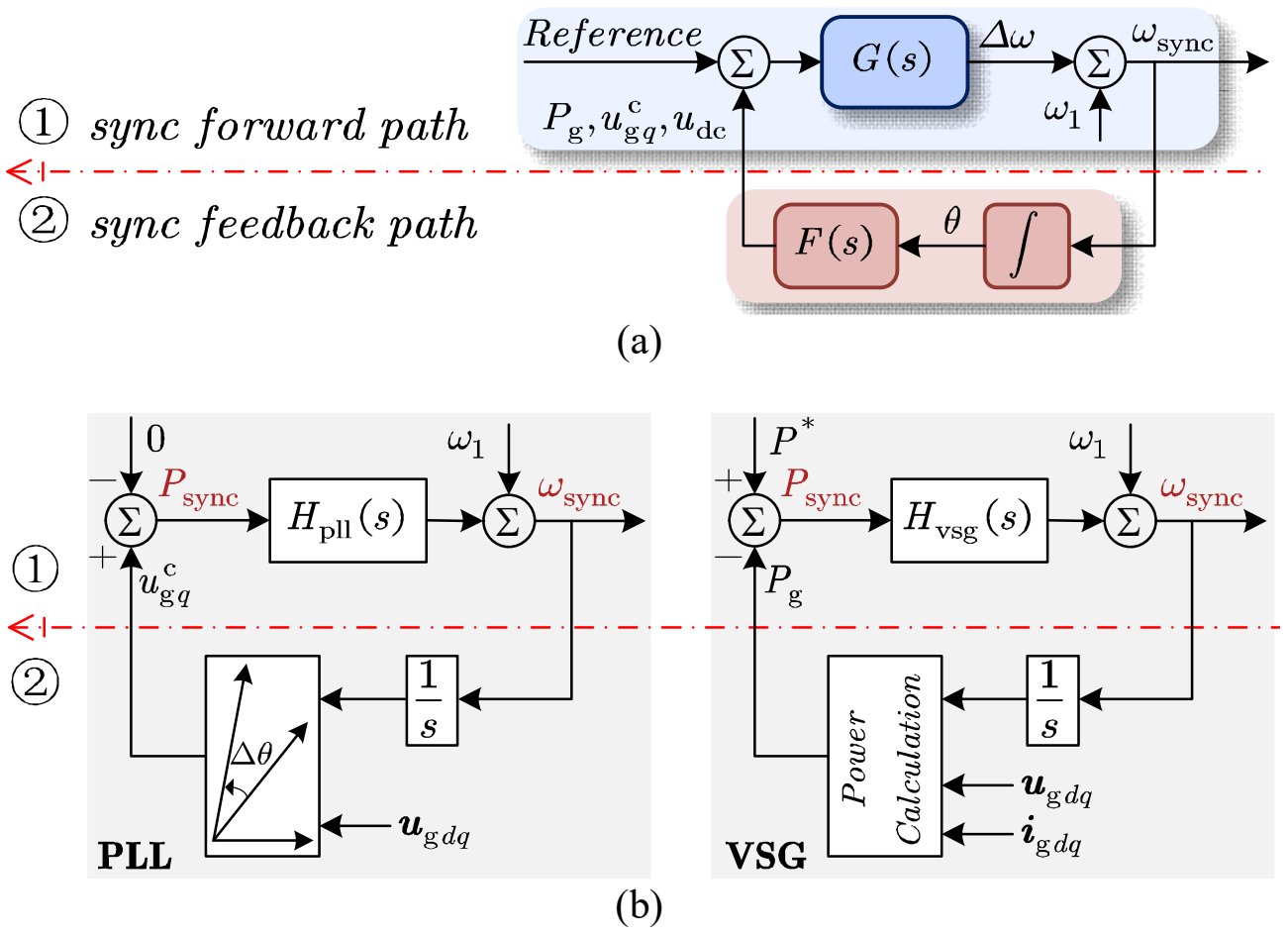

Fig. 3. Unified representation of sync loops, (a) unified feedback form; (b) typical GFL/GFM sync loops.

*2) Derivation of Four-port EIM*

The GFL sync feedback path is influenced by the phase difference between controller and system frame, while the GFM sync feedback path is affected by the active power $P_{\mathrm{g}}$ (can be calculated from voltage and current in ac-side PoC). Their small-signal dynamics can be represented as the function of

$$\begin{aligned} PLL &: \Delta P_{\mathrm{sync}}(s) = G\left(\Delta \boldsymbol{u}_{\mathrm{g}dq}, \Delta\omega_{\mathrm{sync}}\right) \\ VSG &: \Delta P_{\mathrm{sync}}(s) = G\left(\Delta \boldsymbol{u}_{\mathrm{g}dq}, \Delta \boldsymbol{i}_{\mathrm{g}dq}, \Delta\omega_{\mathrm{sync}}\right) \end{aligned} \tag{6}$$

Considering GFM/GFL controls shown in Fig. 4, $\Delta \boldsymbol{i}_{\mathrm{g}dq}$ can be expressed as the function of (details given in Appendix A)

$$\Delta \boldsymbol{i}_{\mathrm{g}dq} = G\left(\Delta \boldsymbol{u}_{\mathrm{g}dq}, \Delta u_{\mathrm{dc}}, \Delta\omega_{\mathrm{sync}}\right) \tag{7}$$

By substituting (7) into the second row of (6), a unified small-signal representation for GFM and GFL sync feedback path can be obtained, by using the same input variables as:

$$\begin{aligned} PLL &: \Delta P_{\mathrm{sync}}(s) = G\left(\Delta \boldsymbol{u}_{\mathrm{g}dq}, 0\cdot\Delta u_{\mathrm{dc}}, \Delta\omega_{\mathrm{sync}}\right) \\ VSG &: \Delta P_{\mathrm{sync}}(s) = G\left(\Delta \boldsymbol{u}_{\mathrm{g}dq}, \Delta u_{\mathrm{dc}}, \Delta\omega_{\mathrm{sync}}\right) \end{aligned} \tag{8}$$

where the input variable $\Delta u_{\mathrm{dc}}$ is also added to the PLL to keep the form unified, but its coefficient is zero. Based on (8), $\Delta P_{\rm sync}$ and $\Delta\omega_{\rm sync}$ of the sync loop will be retained as independent variables, rather than being treated as intermediate variables just like the derivation of the three-port impedance. In this way, four-port EIM in the admittance form can be expressed as:

$$\begin{bmatrix}\Delta\boldsymbol{i}_{gdq}(s)\\ \Delta i_{\rm dc}(s)\\ \Delta P_{\rm sync}(s)\end{bmatrix}=\underbrace{\left[\begin{array}{cc|c}\boldsymbol{Y}_{dq}(s) & \boldsymbol{a}_{dq}(s) & \boldsymbol{c}_{dq}(s)\\ \boldsymbol{b}_{dq}(s) & Y_{\rm dc}(s) & d(s)\\ \hline \boldsymbol{k}_{\rm sync}^{\rm ac}(s) & k_{\rm sync}^{\rm dc}(s) & Y_{\rm sync}^{\rm fe}(s)\end{array}\right]}_{\boldsymbol{Y}_{{\rm con}_dq}^{4}(s)}\begin{bmatrix}\Delta\boldsymbol{u}_{gdq}(s)\\ \Delta u_{\rm dc}(s)\\ \Delta\omega_{\rm sync}(s)\end{bmatrix} \tag{9}$$

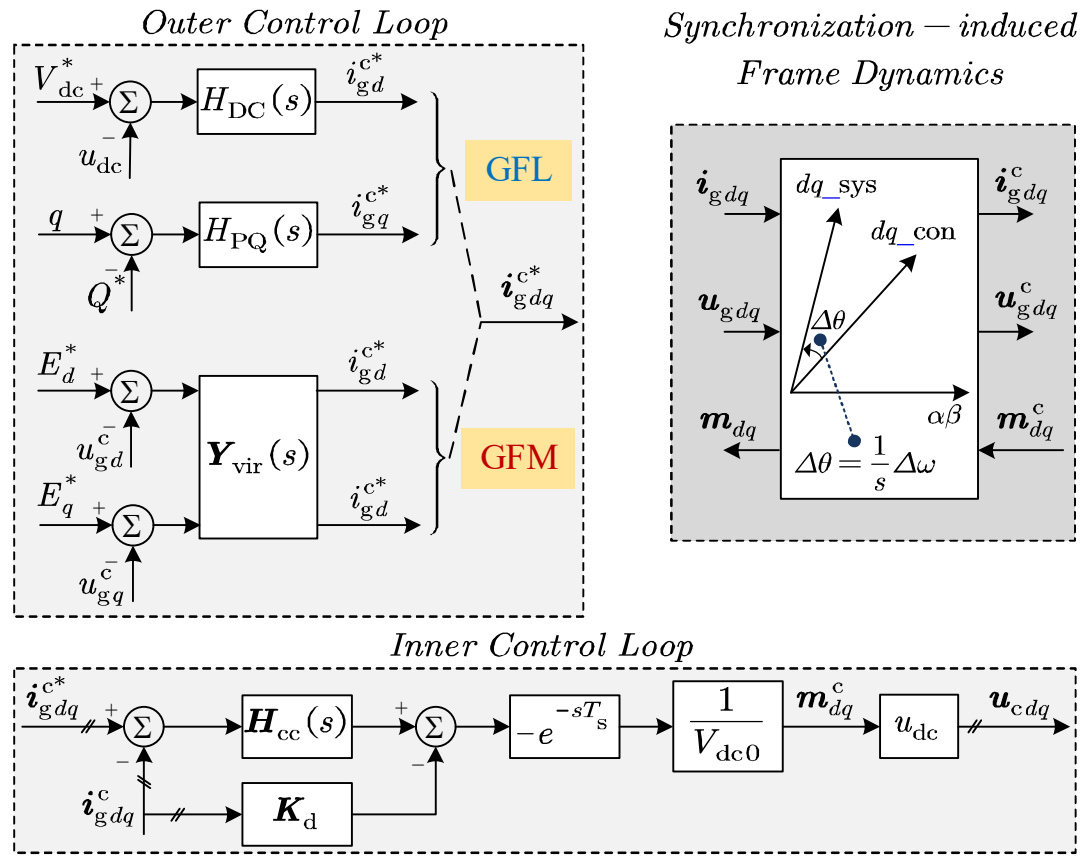

Fig. 4. Diagram of typical GFL/GFM inner and outer control loop.

The upper left block matrix of (9) is exactly in the form of the conventional three-port admittance as (1), but with different analytical expressions due to the separation of sync forward path in (5). A specific derivation of (9) is provided in Appendix A by considering a typical GFL/GFM control given in Fig. 4.

### B. Formulation of the EIN with Sync Nodes

This part will illustrate how to adapt the above EIM with sync port into the formulation of a networked electric system. The resulting model is referred as the EIN as earlier mentioned.

#### 1) Modification of Sync Loop as Electric Branch

An important premise for the network calculus is that the Kirchhoff's voltage/current law (KVL/KCL) should be satisfied [25]. This condition is fulfilled for the ordinary impedance network with only electric nodes/branches, while would fail for separated sync forward path $Z_{\rm sync}^{\rm fo}$ given in (5). To better introduce this property, the couplings of syn loop with $\Delta\boldsymbol{u}_{gdq}$ and $\Delta u_{\rm dc}$ in Fig. 3 are tentatively ignored, and a simplified model is obtained as:

$$\begin{aligned}&\textit{sync forward path}: \Delta\omega_{\rm sync}^{'}=Z_{\rm sync}^{\rm fo}\Delta P_{\rm sync}^{'}\\ &\textit{sync feedback path}: \Delta P_{\rm sync}^{'}=Y_{\rm sync}^{\rm fe}\Delta\omega_{\rm sync}^{'}\end{aligned} \tag{10}$$

where the simplified sync feedback path contains only diagonal element $Y_{\rm sync}^{\rm fe}$ of the four-port EIM in (9). $\Delta P_{\rm sync}^{'}$ and $\Delta\omega_{\rm sync}^{'}$ denote that the influence of $\Delta\boldsymbol{u}_{gdq}$ and $\Delta u_{\rm dc}$ have been excluded from the original $\Delta P_{\rm sync}$ and $\Delta\omega_{\rm sync}$.

The control diagram of (10) is presented in Fig. 5(a). As stated in (5), the input $\Delta P_{\rm sync}^{'}$ and output $\Delta\omega_{\rm sync}^{'}$ of the sync forward path can be regarded as the virtual current and voltage. Based on this, the equivalent circuit of Fig. 5(a) is shown in Fig. 5(b), from which it can be seen that:

- virtual current $\Delta P_{\rm sync}^{'}$ inflows the left-side of the sync forward path $Y_{\rm sync}^{\rm fo}$, and thus the positive direction of virtual voltage $\Delta\omega_{\rm sync}^{'}$ is from left to right;
- virtual current $\Delta P_{\rm sync}^{'}$ outflows the left-side of the sync feedback path $Y_{\rm sync}^{\rm fe}$, and thus the positive direction of virtual voltage $\Delta\omega_{\rm sync}^{'}$ is from right to left.

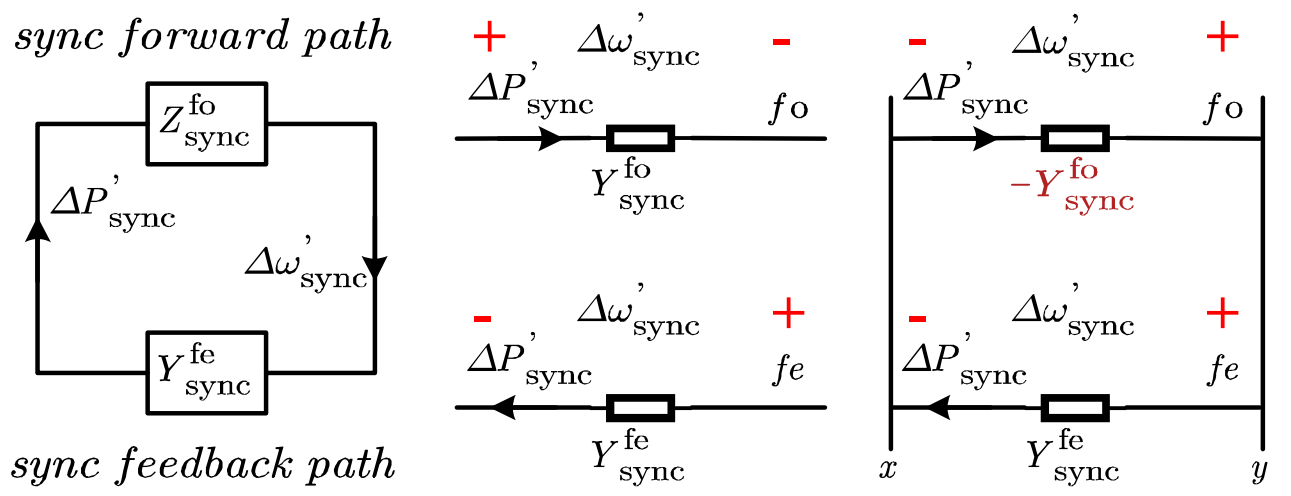

(a) (b) (c)

Fig. 5. Modification of the unified sync loop as directed virtual electric branch, (a) signal flow; (b) polarity of the signal flow in virtual circuit; (c) convention modification of the virtual circuit branch.

Therefore, the positive direction of virtual voltage $\Delta\omega_{\rm sync}^{'}$ for the sync forward path and feedback path is contrary to each other, meaning that these two virtual branches cannot be connected in parallel otherwise the KVL would fail. This is referred as the control convention, and needs to be modified to adhere with the electric convention. As shown in Fig. 5(c), the KVL can be satisfied by modifying $Y_{\rm sync}^{\rm fo}$ as $-Y_{\rm sync}^{\rm fo}$, and then the positive direction of virtual voltage $\Delta\omega_{\rm sync}^{'}$ in sync forward path will be reversed. The above convention modification ensures that the four-port EIM can be interconnected using the circuit law. With this precondition, mature network modeling methods as presented in (2) can be directly applied.

#### 2) Loop Gain Model of the EIN

The EIN modeling will be given first considering a simple case, and then generalized to a generic converters-interlinked system. For the simple case, the equivalent circuit of the single-machine system is shown in Fig. 6(a), where diagonal elements of the four-port EIM constitutes of branch admittances $\boldsymbol{Y}_{dq}$, $Y_{\rm dc}$, $Y_{\rm sync}$, while non-diagonal elements (e.g., $\boldsymbol{a}_{dq}$) represent mutual couplings among branches denoted in black dots. Based on a similar partition concept shown in Fig. 1(a), the EIN of this system can be split into two subsystems at PoCs, where the subsystem consisting of AC, DC and sync nets can be modeled as:

$$\boldsymbol{Z}_{\rm net1}=\left[blkdiag\left(\boldsymbol{Y}_{\rm ac}^{\rm g},Y_{\rm dc}^{\rm g},-Y_{\rm sync}^{\rm fo}\right)\right]^{-1} \tag{11}$$

and the other subsystem consisting of the converter net is in fact the four-port EIM $\boldsymbol{Y}_{{\rm con}_dq}^{4}$ in (9), and is denoted here as $\boldsymbol{Y}_{\rm con1}$.

Based on the above two subsystems' models, the loop gain model of the overall system can be written as:

$$\boldsymbol{L}_1^{\rm EIN}=\boldsymbol{Z}_{\rm net1}\boldsymbol{Y}_{\rm con1}=\begin{bmatrix}\boldsymbol{Z}_{\rm ac}^{\rm g}\boldsymbol{Y}_{dq} & \boldsymbol{Z}_{\rm ac}^{\rm g}\boldsymbol{a}_{dq} & \boldsymbol{Z}_{\rm ac}^{\rm g}\boldsymbol{c}_{dq}\\ Z_{\rm dc}^{\rm g}\boldsymbol{b}_{dq} & Z_{\rm dc}^{\rm g}Y_{\rm dc} & Z_{\rm dc}^{\rm g}d\\ -Z_{\rm sync}^{\rm fo}\boldsymbol{k}_{\rm sync}^{\rm ac} & -Z_{\rm sync}^{\rm fo}k_{\rm sync}^{\rm dc} & -Z_{\rm sync}^{\rm fo}Y_{\rm sync}^{\rm fe}\end{bmatrix} \tag{12}$$

The above modeling procedure can be readily generalized to a multi-converter system in Fig. 6(b). For this generalized system, the subsystem with AC/DC and sync nets can be modeled as:

$$\boldsymbol{Z}_{\rm net}=\begin{bmatrix}\begin{bmatrix}\boldsymbol{Y}_{\rm ac}^{\rm g1} & & \\ & Y_{\rm dc}^{\rm g1} & \\ & & -Y_{\rm sync}^{\rm fo1}\end{bmatrix}_{4\times4} & \cdots & \begin{bmatrix}\boldsymbol{Y}_{\rm ac}^{{\rm g}k1} & & \\ & Y_{\rm dc}^{{\rm g}k1} & \\ & & 0\end{bmatrix}_{4\times4}\\ \vdots & \ddots & \vdots\\ \begin{bmatrix}\boldsymbol{Y}_{\rm ac}^{{\rm g}1k} & & \\ & Y_{\rm dc}^{{\rm g}1k} & \\ & & 0\end{bmatrix}_{4\times4} & \cdots & \begin{bmatrix}\boldsymbol{Y}_{\rm ac}^{{\rm g}k} & & \\ & Y_{\rm dc}^{{\rm g}k} & \\ & & -Y_{\rm sync}^{{\rm fo}k}\end{bmatrix}_{4\times4}\end{bmatrix}^{-1} \tag{13}$$

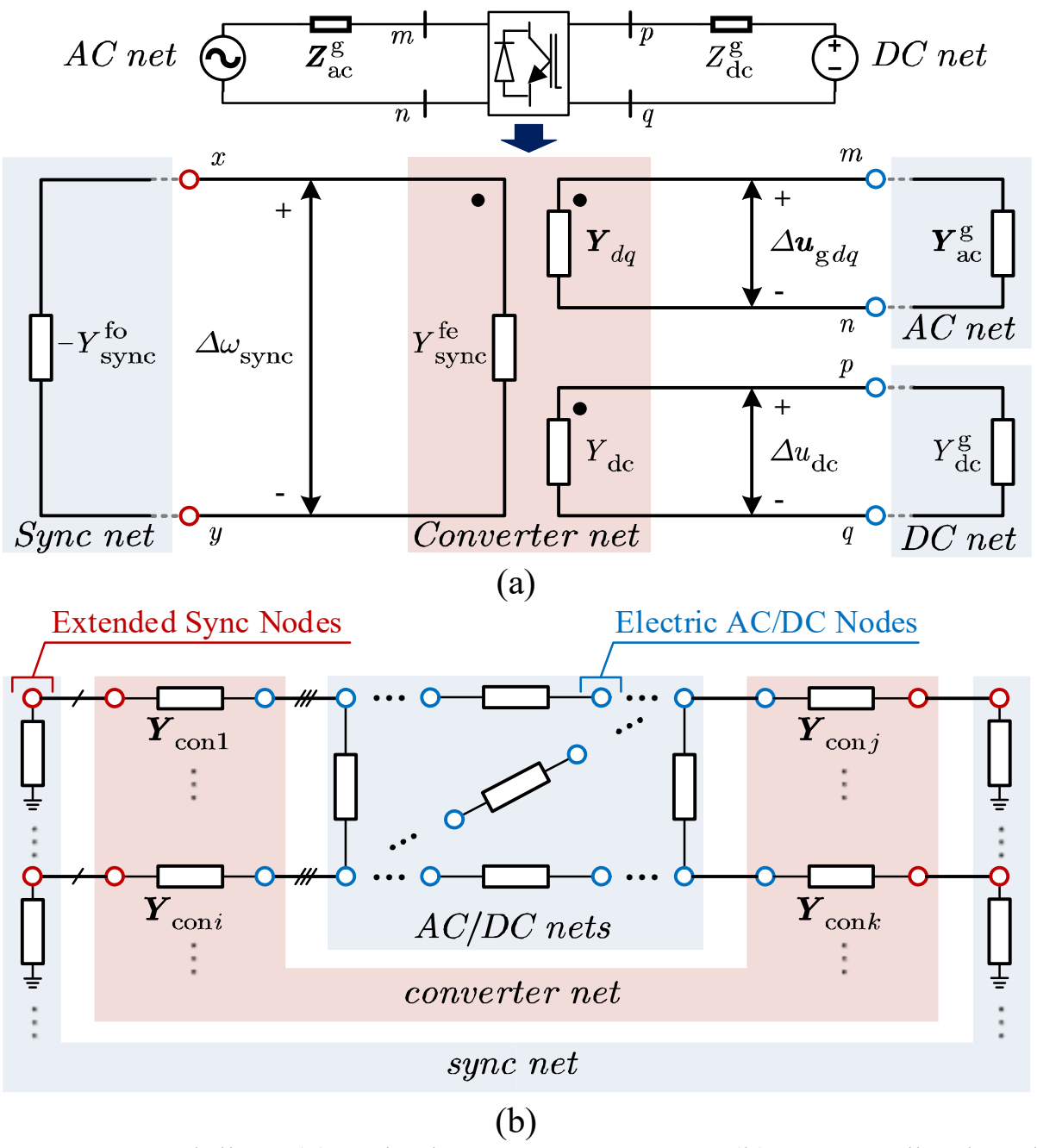

Fig. 6. EIN modeling. (a) a single converter system; (b) a generalized multi-converter system.

where off-diagonal terms like $\boldsymbol{Y}_{ac}^{g1k}$, $Y_{dc}^{g1k}$ are possible mutual admittances among ac and dc nodes, and sync nodes are decoupled from others. The converter subsystem of this case can be modelled as a block diagonal matrix like:

$$\boldsymbol{Y}_{con}=\begin{bmatrix}\left[\boldsymbol{Y}_{con1}\right]_{4\times4} & & \\ & \ddots & \\ & & \left[\boldsymbol{Y}_{conk}\right]_{4\times4}\end{bmatrix} \tag{14}$$

Then, the system loop gain model $\boldsymbol{L}^{\mathrm{EIN}}$ can be obtained by calculating $\boldsymbol{Z}_{net}\,\boldsymbol{Y}_{con}$. It is worth noting that different from the electrical network, the obtained EIN loop gain model contains both electrical and control ports introduced by the four-port EIM, which are of different units. Therefore, a dimensionless transformation is required to convert $\boldsymbol{Z}_{net}$ and $\boldsymbol{Y}_{con}$ of actual unit to those of per unit. The transformation process is given by:

$$per\ unit:\ \boldsymbol{I}_{b}^{-1}\boldsymbol{Y}_{con}\boldsymbol{U}_{b},\ \ \boldsymbol{U}_{b}^{-1}\boldsymbol{Z}_{net}\boldsymbol{I}_{b} \tag{15}$$

where base vector $\boldsymbol{I}_b$ and $\boldsymbol{U}_b$ are defined in Appendix B. The following analysis is all based on the per-unit model.

### C. Sync-Node Extended Sensitivity Analysis

By substituting the above EIN model $\boldsymbol{L}^{\mathrm{EIN}}$ into (3) and (4), the sync-node extended PF matrix can be obtained. Differed with the PF of the ordinary impedance network, sync dynamics are incorporated into such PF matrix, typically the sensitivity of the diagonal element $Z_{sync}^{fo}\,Y_{sync}^{fe}$ to the studied mode $\lambda_m$. Based on the derivative chain rule, the sensitivity of the interested mode $\lambda_m$ with respect to each component inside $\boldsymbol{L}^{\mathrm{EIN}}$ can be further extracted, including sync virtual branch $Z_{sync}^{fo}$ and $Y_{sync}^{fe}$. As to components inside $\boldsymbol{Z}_{net}$ (e.g., sync forward virtual branch $Z_{sync}^{fo}$), their sensitivities can be calculated by (16), where $Z_{ij}$ in $i$ row and $j$ column of $\boldsymbol{Z}_{net}$ is used for the concrete derivation.

$$\begin{aligned}S_{ij}^{Z}&=\frac{\partial \Lambda_k}{\partial Z_{ij}}=\boldsymbol{t}_k\frac{\partial \boldsymbol{L}^{\mathrm{EIN}}}{\partial Z_{ij}}\boldsymbol{r}_k\\&=\boldsymbol{t}_k\left(\frac{\partial \boldsymbol{Z}_{net}}{\partial Z_{ij}}\boldsymbol{Y}_{con}+\boldsymbol{Z}_{net}\frac{\partial \boldsymbol{Y}_{con}}{\partial Z_{ij}}\right)\boldsymbol{r}_k\end{aligned} \tag{16}$$

Since $\boldsymbol{Y}_{con}$ do not contain any element concerning $Z_{ij}$, i.e., partial derivatives equal to zero, (16) can be simplified as:

$$\begin{aligned}S_{ij}^{Z}&=\boldsymbol{t}_k\left(\frac{\partial \boldsymbol{Z}_{net}}{\partial Z_{ij}}\boldsymbol{Y}_{con}\right)\boldsymbol{r}_k=\boldsymbol{t}_k\left(\begin{bmatrix}\ddots & 0 & \iddots\\ 0 & 1_{ij} & 0\\ \iddots & 0 & \ddots\end{bmatrix}\cdot\boldsymbol{Y}_{con}\right)\boldsymbol{r}_k\\&=\boldsymbol{t}_k\left(\begin{matrix}\vdots\\ i\\ \vdots\end{matrix}\begin{bmatrix}\ddots & 0 & \iddots\\ Y_{j1} & \cdots & Y_{jn}\\ \iddots & 0 & \ddots\end{bmatrix}\right)\boldsymbol{r}_k=\begin{bmatrix}t_{ki}Y_{j1} & \cdots & t_{ki}Y_{jn}\end{bmatrix}\boldsymbol{r}_k\\&=t_{ki}Y_{j1}r_{1k}+\cdots+t_{ki}Y_{jn}r_{nk}=\left[\boldsymbol{Y}_{con}\cdot\boldsymbol{PF}^{L}\right]_{ji}\end{aligned} \tag{17}$$

For components inside $\boldsymbol{Y}_{con}$ (e.g., sync feedback virtual branch $Y_{sync}^{fe}$), their sensitivities can be calculated by (18), where $Y_{ij}$ in $i$ row and $j$ column of $\boldsymbol{Y}_{con}$ is taken for example.

$$\begin{aligned}S_{ij}^{Y}&=\boldsymbol{t}_k\left(\boldsymbol{Z}_{net}\frac{\partial \boldsymbol{Y}_{con}}{\partial Y_{ij}}\right)\boldsymbol{r}_k=\boldsymbol{t}_k\left(\boldsymbol{Z}_{net}\cdot\begin{bmatrix}\ddots & 0 & \iddots\\ 0 & 1_{ij} & 0\\ \iddots & 0 & \ddots\end{bmatrix}\right)\boldsymbol{r}_k\\&=t_{k1}Z_{1i}r_{jk}+\cdots+t_{kn}Z_{ni}r_{jk}=\left[\boldsymbol{PF}\cdot\boldsymbol{Z}_{net}\right]_{ji}\end{aligned} \tag{18}$$

## IV. Case Study and Validation

### A. Test System Description

The test hybrid GFL/GFM HVDC system is shown in Fig. 7(a), where sending converter (SEC) operates with GFM mode and receiving converter (REC) operates with GFL mode. The EIN of the studied system is given in Fig. 7(b), where node #1, #6 are virtual sync nodes, and the rest are electrical nodes. The impedance $\boldsymbol{Z}_{net}(s)$ of the subsystem containing AC, DC and sync nets can be modeled as:

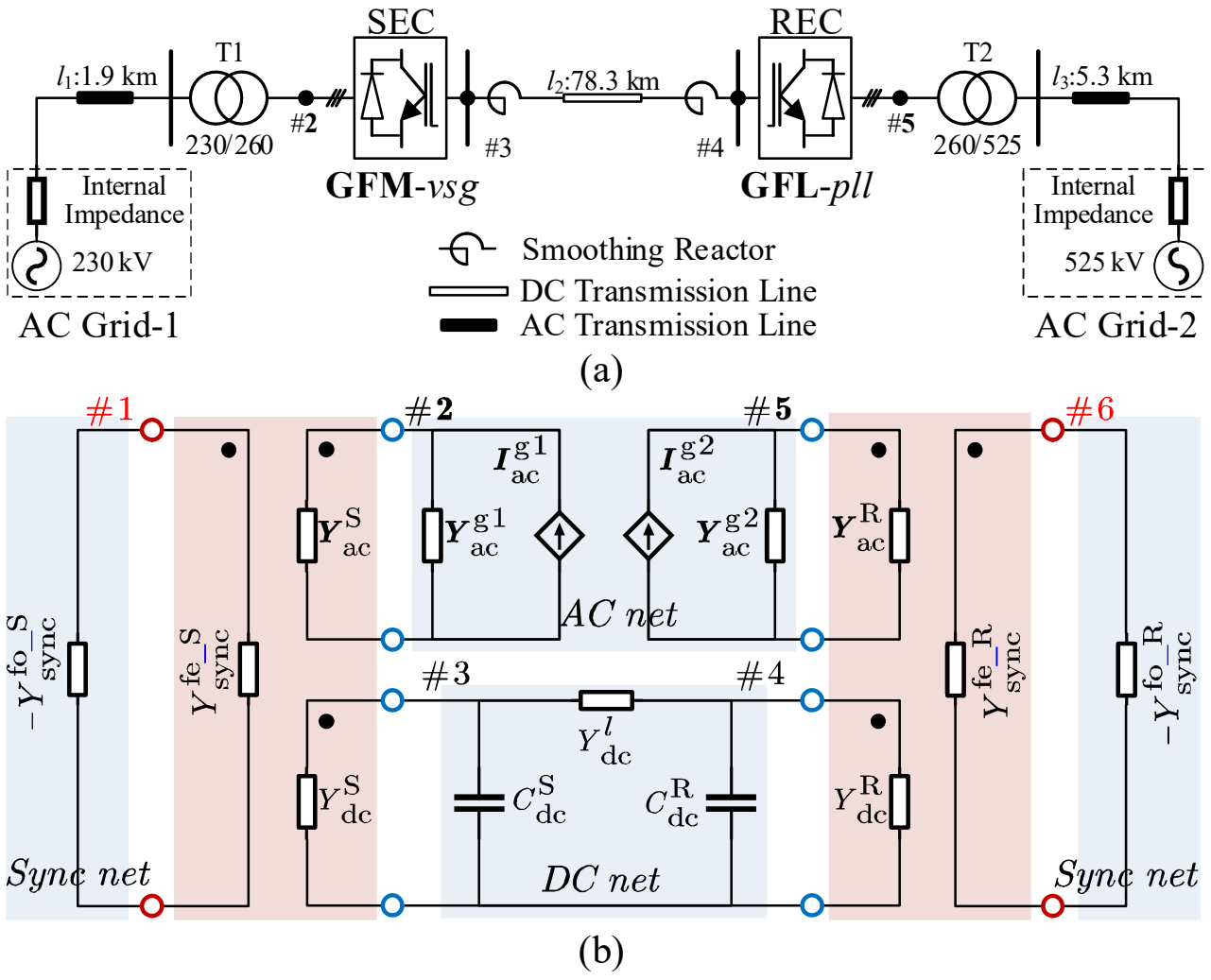

Fig. 7. Test system, (a) system topology; (b) equivalent circuit with modified sync virtual branches.

$$Z_{\text{net}}(s)=\left[blkdiag\left(-Y_{\text{sync}}^{\text{fo_S}},\boldsymbol{Y}_{\text{ac}}^{\text{g1}},\boldsymbol{Y}_{\text{dc}}^{\text{g}},\boldsymbol{Y}_{\text{ac}}^{\text{g2}},-Y_{\text{sync}}^{\text{fo_R}}\right)\right]^{-1} \quad (19)$$

where $\boldsymbol{Y}_{\text{dc}}^{\text{g}}=\begin{bmatrix} sC_{\text{dc}}^{\text{S}}+Y_{\text{dc}}^{l} & -Y_{\text{dc}}^{l} \\ -Y_{\text{dc}}^{l} & sC_{\text{dc}}^{\text{R}}+Y_{\text{dc}}^{l}\end{bmatrix}$.

The admittance $\boldsymbol{Y}_{\text{con}}(s)$ of the subsystem containing SEC and REC can be represented as in (20). Then, the system loop gain model can be obtained as $\boldsymbol{L}(s)=\boldsymbol{Z}_{\text{net}}(s)\,\boldsymbol{Y}_{\text{con}}(s)$.

$$\boldsymbol{Y}_{\text{con}}(s)=blkdiag\left\{\boldsymbol{Y}_{\text{con}}^{\text{SEC}},\boldsymbol{Y}_{\text{con}}^{\text{REC}}\right\}=$$
$$blkdiag\left(\begin{bmatrix} Y_{\text{sync}}^{\text{fe_S}} & \boldsymbol{k}_{\text{sync}}^{\text{ac_S}} & k_{\text{sync}}^{\text{dc_S}} \\ \boldsymbol{c}^{\text{S}} & \boldsymbol{Y}_{\text{ac}}^{\text{S}} & \boldsymbol{a}^{\text{S}} \\ d^{\text{S}} & \boldsymbol{b}^{\text{S}} & Y_{\text{dc}}^{\text{S}}\end{bmatrix},\begin{bmatrix} Y_{\text{dc}}^{\text{R}} & \boldsymbol{b}^{\text{R}} & d^{\text{R}} \\ \boldsymbol{a}^{\text{R}} & \boldsymbol{Y}_{\text{ac}}^{\text{R}} & \boldsymbol{c}^{\text{R}} \\ k_{\text{sync}}^{\text{dc_R}} & \boldsymbol{k}_{\text{sync}}^{\text{ac_R}} & Y_{\text{sync}}^{\text{fe_R}}\end{bmatrix}\right) \quad (20)$$

## B. *Validation of the Four-port EIM*

The accuracy of the four-port EIM will be validated by frequency scanning results. To acquire the measured four-port EIM, two type of VSCs in the test system need to be connected with ideal sources as shown in Fig. 8, so as to acquire their open-loop impedances. Then, totally four independent perturbations are injected from the main circuit and control part, two of which are from the ac-side and the third one is from the dc-side. This is identical to the measurement of the classical three-port impedance, while the difference lies in the fourth perturbation, in which an additional injection from the sync control port is applied and corresponding response (i.e., $\omega_{\text{sync}}$) is collected. In this process, the applied injections (e.g., series voltage) are in the form of:

$$v_{inj1,2}=V_{inj}^{\text{ac}}\begin{bmatrix}\sin\left(\left[\omega_{inj}\pm\omega_1\right]t\right)\\ \sin\left(\left[\omega_{inj}\pm\omega_1\right]t\mp\frac{2}{3}\pi\right)\\ \sin\left(\left[\omega_{inj}\pm\omega_1\right]t\pm\frac{2}{3}\pi\right)\end{bmatrix}$$
$$v_{inj3}=V_{inj}^{\text{dc}}\sin\left(\omega_{inj}t\right) \quad (21)$$
$$v_{inj4}=P_{inj}^{\text{sync}}\sin\left(\omega_{inj}t\right)$$

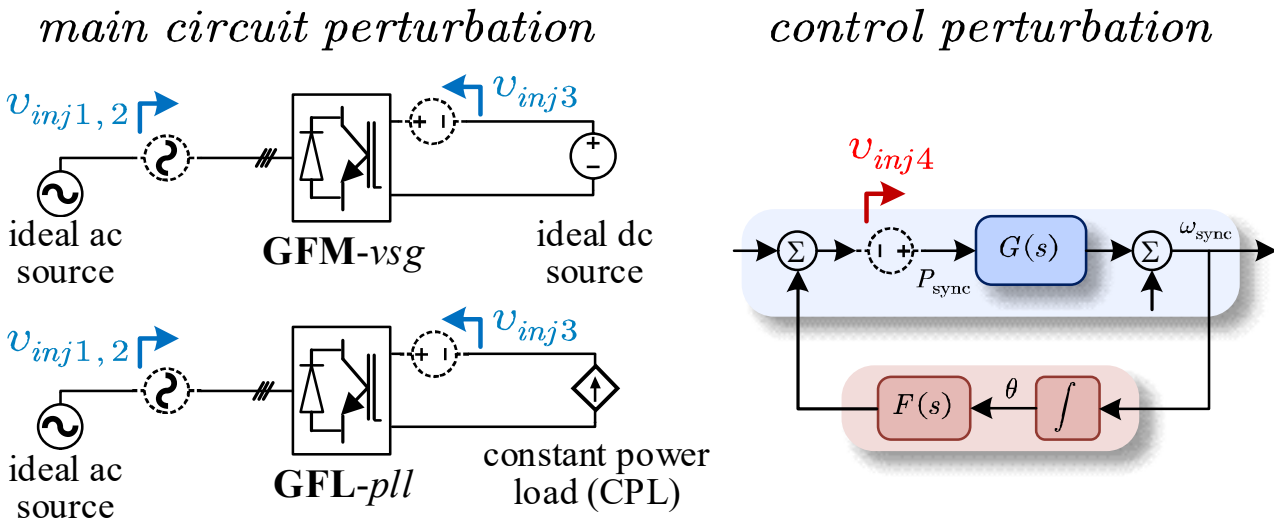

Fig. 8. Perturbations applied in the frequency scanning method.

The measured four-port EIM can be calculated as:

$$\boldsymbol{Y}_{\text{con_}pn}^{4}=\begin{bmatrix} I_{p1} & I_{p2} & I_{p3} & I_{p4} \\ I_{n1} & I_{n2} & I_{n3} & I_{n4} \\ I_{\text{dc}1} & I_{\text{dc}2} & I_{\text{dc}3} & I_{\text{dc}4} \\ P_1 & P_2 & P_3 & P_4\end{bmatrix}\begin{bmatrix} V_{p1} & V_{p2} & V_{p3} & V_{p4} \\ V_{n1} & V_{n2} & V_{n3} & V_{n4} \\ V_{\text{dc}1} & V_{\text{dc}2} & V_{\text{dc}3} & V_{\text{dc}4} \\ \omega_1 & \omega_2 & \omega_3 & \omega_4\end{bmatrix}^{-1} \quad (22)$$

where (22) is in modified sequence domain [20], which can be obtained from *dq* domain by applying (23).

$$\boldsymbol{Y}_{\text{con_}pn}^{4}=\boldsymbol{A}_Z\cdot\boldsymbol{Y}_{\text{con_}dq}^{4}\cdot\boldsymbol{A}_Z^{-1} \quad (23)$$

where $\boldsymbol{A}_Z=\begin{bmatrix}\boldsymbol{a}_z & \\ & \boldsymbol{I}\end{bmatrix}$, $\boldsymbol{a}_z=\frac{1}{\sqrt{2}}\begin{bmatrix}1 & j\\ 1 & -j\end{bmatrix}$, and $\boldsymbol{I}=\begin{bmatrix}1 & 0\\ 0 & 1\end{bmatrix}$. The comparison results are given in Fig. 9(a) and (b), including GFL PLL and GFM VSG. Theoretical curves match very well with measurement curves, which verify the accuracy of the proposed four-port EIM.

## C. *Verification of Sync Node-Extended Modal Analysis*

### 1) *Stability Judgement*

Three representative oscillation cases are analyzed, including low sync damping coefficient, high PLL bandwidth and weak grid, whose information are given in Fig. 10(a)~(b). It can be seen that for all three cases, theoretical results conform well with time-domain simulations, and thus validate the accuracy of the established EIN model. Taking the Case I for an intuitive illustration, by decreasing the damping coefficient *D*, the EIN model will generate one pair of RHP-zeros at around 4.3974Hz (*dq* frame), indicating that the system is unstable. This is consistent with simulation results, where the system went from being stable to unstable for the same marginal condition and the oscillation frequency is about 4.3333 Hz (*dq* frame).

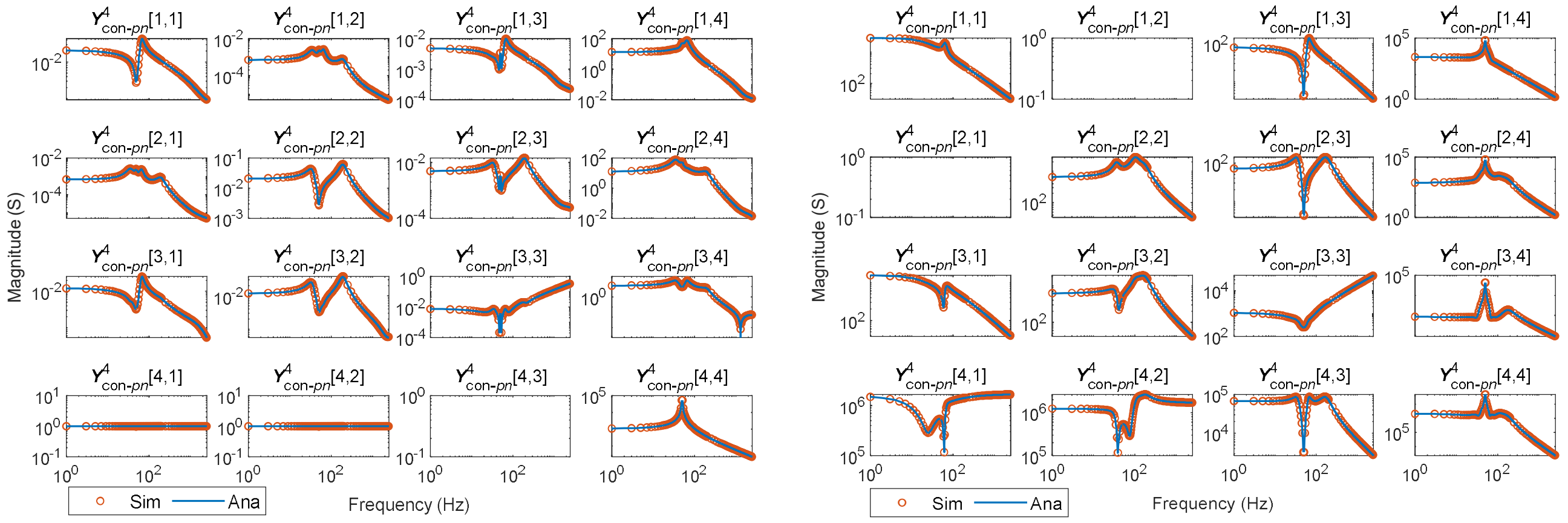

(a) (b)

Fig. 9. Frequency scan-based validation of the proposed Four-port EIM, (a) GFL with PLL; (b) GFM with VSG.

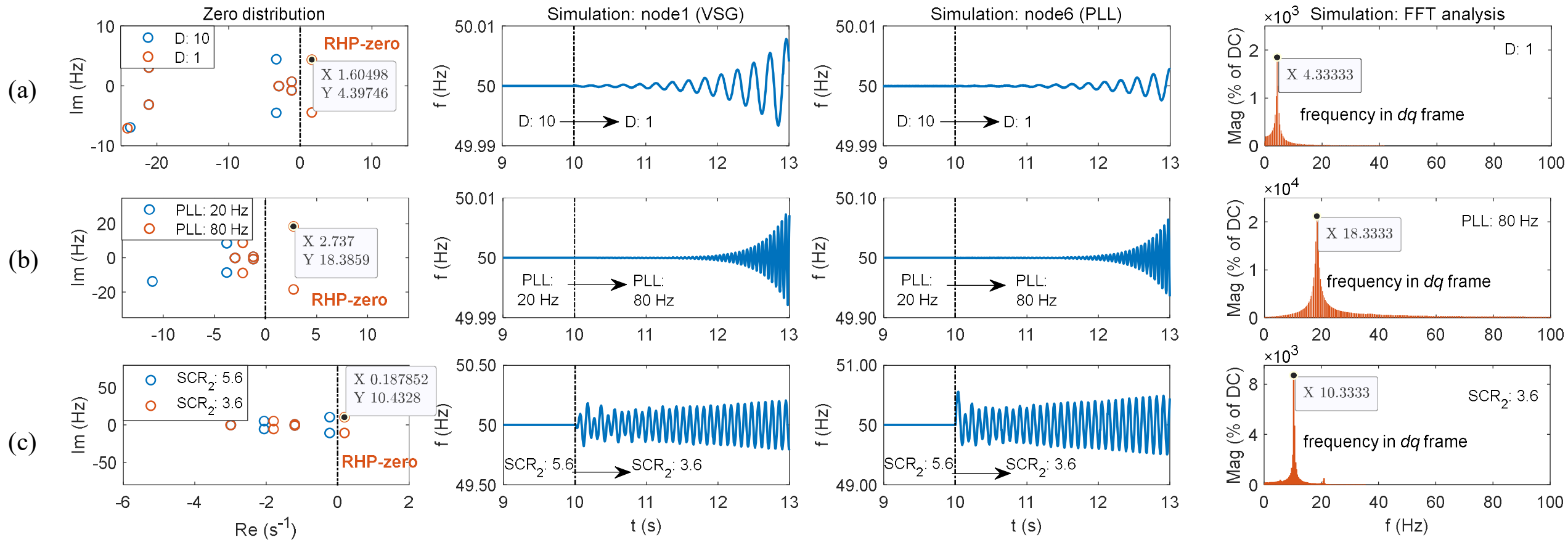

Fig. 10. Oscillation cases, (a) Case I: reduce VSG damping coefficient in SEC from 10 p.u. to 1 p.u.; (b) Case II: increase PLL control bandwidth in REC from 20Hz to 80Hz; (c) Case III: reduce the SCR of AC grid-2 from 5.6 to 3.6.

### 2) *Modal Sensitivity Analysis*

The node PF of Case I, II and III is calculated in Table I. And the sensitivity of each impedance and admittance component shown in (19) and (20) is calculated in Fig. 11. Taking Case I for instance, the accuracy of the component sensitivity is also validated in Table II. It can be seen that under a 5% increment in each component, errors between the predicted and actual $\Delta\Lambda_k$ are within 5%. The predicted change $\Delta\Lambda_k$ is calculated by multiplying the theoretically calculated sensitivity with the corresponding component 5% increment, while the actual change is obtained by re-computing $\Lambda_k$ of the updated system model (replace the corresponding component value with 105% times of original value) and subtracting the original $\Lambda_k$.

It can be seen from Table I and Fig. 11 that influences of sync dynamics originally implicit inside converters on oscillations can now be intuitively analyzed via the extended sync nodes (node #1 and #6) and virtual branches ($Z_{\rm sync}^{\rm fo_S}$, $Z_{\rm sync}^{\rm fo_R}$, $Y_{\rm sync}^{\rm fe_S}$ and $Y_{\rm sync}^{\rm fe_R}$), just as other electrical nodes and branches. For Case I, virtual sync node #1 and electrical node #2 are with highest PFs as shown in Table I. Furthermore, it can be obtained from Fig. 11(a) that the sync forward virtual branch $Z_{\rm sync}^{\rm fo_S}$ of VSG has most significant impact on the system oscillation, indicating that the system instability is caused by VSG sync control.

For Case II, electrical node #5 and virtual sync node #6 are with highest PFs as shown in Table I. And according to Fig. 11(b), the sync forward virtual branch $Z_{\rm sync}^{\rm fo_R}$ of PLL and AC grid-2 impedance $\boldsymbol{Z}_{\rm ac}^{\rm g2}$ will impose great influence on the system oscillation, which means that the system instability is resulted from the interaction between PLL and weak ac grid. This conclusion is consistent with the current recognition of the weak grid-induced instability in GFL system [6].

For Case III, the oscillation is triggered by changing the AC grid-2 SCR, and electrical node #5 is with highest PF as shown in Table I. It is further observed from Fig. 11(c) that VSG sync forward virtual branch $Z_{\rm sync}^{\rm fo_S}$, PLL sync feedback virtual branch $Z_{\rm sync}^{\rm fo_R}$, AC grid-1 impedance $\boldsymbol{Z}_{\rm ac}^{\rm g1}$ and AC grid-2 impedance $\boldsymbol{Z}_{\rm ac}^{\rm g2}$ are with high component sensitivity. This indicates that the system instability is caused by the complex interaction among PLL, VSG and ac grids, which can be intuitively reflected by the sync-node extended sensitivity indices.

The above three typical cases intuitively demonstrate the effectiveness of the sync node extended FMA in analyzing sync dynamics and their interactions on system oscillations.

TABLE I
PARTICIPATION FACTOR

| Node | PF (mag.) | Node | PF (mag.) | Node | PF (mag.) |
|---|---|---|---|---|---|
| **Case I** | | | | | |
| **#1** | **0.7181** | #3 | 0.0020 | #5 | 0.0004 |
| **#2** | **0.2779** | #4 | 0.0015 | #6 | 0.0001 |
| **Case II** | | | | | |
| 1 | 0.0001 | 3 | 0.0006 | **5** | **0.6218** |
| 2 | 0.0035 | 4 | 0.1163 | **6** | **0.2577** |
| **Case III** | | | | | |
| 1 | 0.0041 | 3 | 0.0368 | **5** | **0.7747** |
| 2 | 0.0511 | 4 | 0.0793 | 6 | 0.0541 |

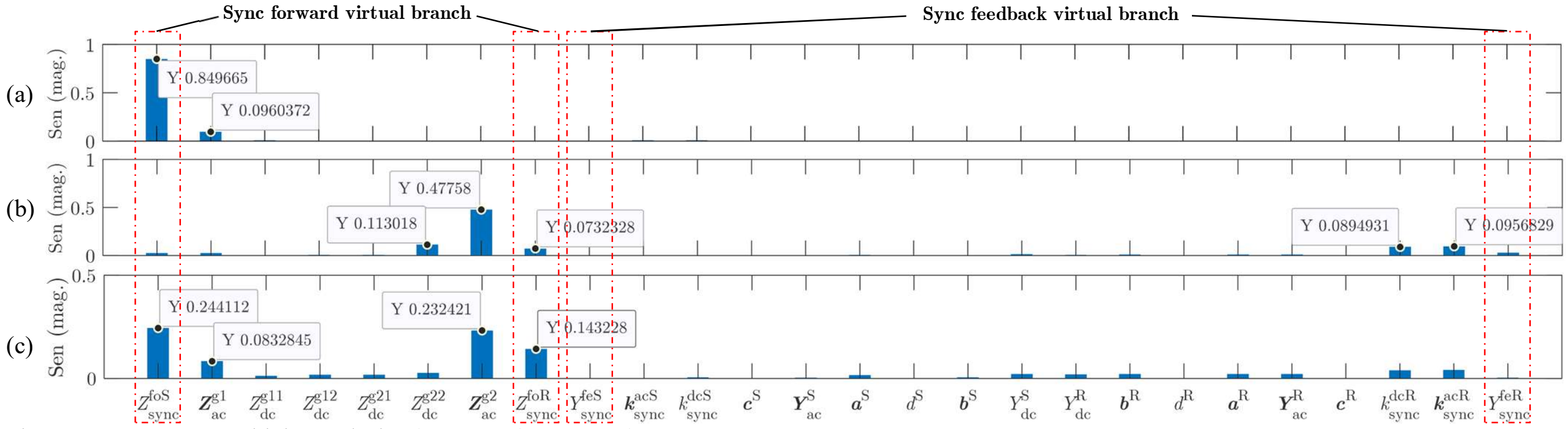

Fig. 11. Component sensitivity analysis, 1) Case I; 2) Case II; 3) Case III.

TABLE II
COMPONENT SENSITIVITY AND TUNING RESULTS UNDER 5% INCREMENT

| Component | Sensitivity | Predicted $\lvert\Delta\Lambda'_m\rvert$ | Actual $\lvert\Delta\Lambda'_m\rvert$ | Error |
|---|---|---|---|---|
| $Z_{\text{sync}}^{\text{fo_S}}$ | 0.2226+$j$2.1442 | 0.1280 | 0.1299 | 1.47% |
| $\boldsymbol{Z}_{\text{ac}}^{\text{g1_11}}$ | 0.1382+$j$0.0381 | 0.0082 | 0.0083 | 1.46% |
| $\boldsymbol{Z}_{\text{ac}}^{\text{g1_12}}$ | 0.0351-$j$0.0211 | 0.1373 | 0.1368 | 0.36% |
| $\boldsymbol{Z}_{\text{ac}}^{\text{g1_21}}$ | -0.1841+$j$0.0209 | 0.6221 | 0.6194 | 0.43% |
| $\boldsymbol{Z}_{\text{ac}}^{\text{g1_22}}$ | -0.0319+$j$0.0422 | 0.0030 | 0.0031 | 1.48% |
| $\boldsymbol{Z}_{\text{dc}}^{\text{g_11}}$ | 0.0157+$j$0.0000 | 0.0206 | 0.0209 | 1.18% |
| $\boldsymbol{Z}_{\text{dc}}^{\text{g_21}}$ | 0.0123+$j$0.0000 | 0.0150 | 0.0152 | 1.68% |
| $\boldsymbol{Z}_{\text{dc}}^{\text{g_12}}$ | 0.0123+$j$0.0000 | 0.0150 | 0.0152 | 1.68% |
| $\boldsymbol{Z}_{\text{dc}}^{\text{g_22}}$ | 0.0096+$j$0.0000 | 0.0126 | 0.0128 | 1.32% |
| $\boldsymbol{Z}_{\text{ac}}^{\text{g2_11}}$ | 0.0061-$j$0.0081 | 0.0028 | 0.0028 | 1.43% |
| $\boldsymbol{Z}_{\text{ac}}^{\text{g2_12}}$ | -0.0004-$j$0.0048 | 0.0122 | 0.0124 | 1.78% |
| $\boldsymbol{Z}_{\text{ac}}^{\text{g2_21}}$ | -0.0016+$j$0.0048 | 0.0127 | 0.0130 | 1.68% |
| $\boldsymbol{Z}_{\text{ac}}^{\text{g2_22}}$ | 0.0009+$j$0.0022 | 0.0007 | 0.0007 | 1.67% |
| $Z_{\text{sync}}^{\text{fo_R}}$ | -0.0010+$j$0.0041 | 0.0098 | 0.0095 | 3.55% |

## V. CONCLUSION

This paper proposed a sync node-extended FMA method for better studying participations of sync dynamics in oscillations of GFL/GFM converters-interlinked systems. In this method, a four-port EIM of the converter is achieved by explicitly modeling the sync loop as a new port and further augmented to the typical three-port impedance model. Then, the sync loop is virtualized as a circuit branch such that the four-port EIM can be adapted to the convention of system-level circuit modeling. Given by this preserved scalability, the EIN of a generic multi-converter system is demonstrated to be readily constructed by networked circuit modeling (e.g., using the node admittance method). The proposed four port model is verified by frequency scanning, and the validity of the EIN-based modal analysis in better identifying the participation of GFL/GFM sync loops in oscillation studies is certified by time domain simulations.

The presented sync node-extended modeling method is an extension to the existing impedance modeling framework. As shown in this work, its conjunction with the FMA method forms a promising tool for oscillation diagnosis, particularly its capacity of quantifying the impacts of sync dynamics.

## APPENDIX A
## CONCRETE DERIVATION OF THE FOUR-PORT EIM

### A. Modeling of Main Circuit and Frame Transformation

The main circuit of VSC can be modeled as:

$$\Delta\boldsymbol{u}_{\text{c}dq}(s) = -\boldsymbol{Z}_{dq}^{\text{f}}(s)\Delta\boldsymbol{i}_{\text{g}dq}(s) + \Delta\boldsymbol{u}_{\text{g}dq}(s) \tag{24}$$

where $\boldsymbol{Z}_{dq}^{\text{f}}(s)$ denotes the impedance of AC filter. VSC's ac-side is coupled with dc-side via the AC/DC power balance, i.e., $\Delta P_{\text{c}}^{\text{dc}}(s)=\Delta P_{\text{c}}^{\text{ac}}(s)$, which can be modeled as:

$$\Delta i_{\text{dc}}(s) = \underbrace{\frac{3\boldsymbol{U}_{\text{c}dq0}}{2V_{\text{dc0}}}}_{\boldsymbol{k}_{\text{dc1}}}\Delta\boldsymbol{i}_{\text{g}dq}(s) + \underbrace{\frac{3\boldsymbol{I}_{\text{g}dq0}}{2V_{\text{dc0}}}}_{\boldsymbol{k}_{\text{dc2}}}\Delta\boldsymbol{u}_{\text{c}dq}(s) + \underbrace{-\frac{I_{\text{dc0}}}{V_{\text{dc0}}}}_{k_{\text{dc3}}}\Delta u_{\text{dc}}(s) \tag{25}$$

where $\boldsymbol{U}_{\text{c}dq0} = \begin{bmatrix} U_{\text{c}d0} & U_{\text{c}q0} \end{bmatrix}$ and $\boldsymbol{I}_{\text{g}dq0} = \begin{bmatrix} I_{\text{g}d0} & I_{\text{g}q0} \end{bmatrix}$. The relationship between the controller frame and the system frame can be represented as:

$$\begin{aligned} \Delta\boldsymbol{i}_{\text{g}dq}^{\text{c}}(s) &= \Delta\boldsymbol{i}_{\text{g}dq}(s) + \boldsymbol{I}_{\text{g}dq0}^{\text{sync}}(s)\Delta\omega_{\text{sync}}(s) \\ \Delta\boldsymbol{u}_{\text{g}dq}^{\text{c}}(s) &= \Delta\boldsymbol{u}_{\text{g}dq}(s) + \boldsymbol{U}_{\text{g}dq0}^{\text{sync}}(s)\Delta\omega_{\text{sync}}(s) \\ \Delta\boldsymbol{m}_{dq}^{\text{c}}(s) &= \Delta\boldsymbol{m}_{dq}(s) + \boldsymbol{M}_{dq0}^{\text{sync}}(s)\Delta\omega_{\text{sync}}(s) \end{aligned} \tag{26}$$

where variables with superscript 'c' are in the controller frame, and variables without superscripts are in the system frame. And

$$\boldsymbol{I}_{\text{g}dq0}^{\text{sync}}(s) = \frac{1}{s}\begin{bmatrix} I_{\text{g}q0} \\ -I_{\text{g}d0} \end{bmatrix},\ \boldsymbol{U}_{\text{g}dq0}^{\text{sync}}(s) = \frac{1}{s}\begin{bmatrix} U_{\text{g}q0} \\ -U_{\text{g}d0} \end{bmatrix},\ \boldsymbol{M}_{dq0}^{\text{sync}}(s) = \frac{1}{s}\begin{bmatrix} m_{q0} \\ -m_{d0} \end{bmatrix}.$$

For PLL, $\Delta P_{\text{sync}}=\Delta u_{\text{g}q}^{\text{c}}$ and can be modeled as in (27). And for VSG , $\Delta P_{\text{sync}}=-\Delta P_{\text{g}}$ and can be modeled as in (28).

$$\Delta P_{\text{sync}}(s) = \boldsymbol{I}_1\Delta\boldsymbol{u}_{\text{g}dq}(s) + \boldsymbol{I}_1\boldsymbol{U}_{\text{g}dq0}^{\text{sync}}(s)\Delta\omega_{\text{sync}}(s) \tag{27}$$

$$\Delta P_{\text{sync}}(s) = -\boldsymbol{I}_2\boldsymbol{U}_{\text{g}dq0}^{\text{PQ}}\Delta\boldsymbol{i}_{\text{g}dq}(s) - \boldsymbol{I}_2\boldsymbol{I}_{\text{g}dq0}^{\text{PQ}}\Delta\boldsymbol{u}_{\text{g}dq}(s) \tag{28}$$

where $\boldsymbol{U}_{\text{g}dq0}^{\text{PQ}} = \frac{3}{2}\begin{bmatrix} U_{\text{g}d0} & U_{\text{g}q0} \\ U_{\text{g}q0} & -U_{\text{g}d0} \end{bmatrix}$, $\boldsymbol{I}_{\text{g}dq0}^{\text{PQ}} = \frac{3}{2}\begin{bmatrix} I_{\text{g}d0} & I_{\text{g}q0} \\ -I_{\text{g}q0} & I_{\text{g}d0} \end{bmatrix}$, $\boldsymbol{I}_1 = \begin{bmatrix} 0 & 1 \end{bmatrix}$ and $\boldsymbol{I}_2 = \begin{bmatrix} 1 & 0 \end{bmatrix}$.

### B. Modeling of Control Part

For GFL outer loop, the output $\Delta\boldsymbol{i}_{\text{g}dq}^{\text{c}*}(s)$ can be modeled as:

$$\begin{aligned} \Delta\boldsymbol{i}_{\text{g}dq}^{\text{c}*}(s) &= \boldsymbol{H}_{\text{dc}}(s)\Delta u_{\text{dc}}(s) + \boldsymbol{H}_{\text{PQ}}(s)\Delta\boldsymbol{PQ}_{\text{g}} = \\ &\underbrace{\boldsymbol{H}_{\text{PQ}}\boldsymbol{U}_{\text{g}dq0}^{\text{PQ}}}_{\boldsymbol{T}_{\text{DC1}}(s)}\Delta\boldsymbol{i}_{\text{g}dq}(s) + \underbrace{\boldsymbol{H}_{\text{PQ}}\boldsymbol{I}_{\text{g}dq0}^{\text{PQ}}}_{\boldsymbol{T}_{\text{DC2}}(s)}\Delta\boldsymbol{u}_{\text{g}dq}(s) + \boldsymbol{H}_{\text{dc}}(s)\Delta u_{\text{dc}}(s) \end{aligned} \tag{29}$$

where $\boldsymbol{H}_{\text{dc}}(s) = \begin{bmatrix} -H_{\text{dc}}(s) \\ 0 \end{bmatrix}$, $\boldsymbol{H}_{\text{PQ}}(s) = \begin{bmatrix} 0 & \\ & H_{\text{PQ}}(s) \end{bmatrix}$.

For GFM outer loop, the output $\Delta\boldsymbol{i}_{\text{g}dq}^{\text{c}*}$ can be modeled as:

$$\Delta\boldsymbol{i}_{\text{g}dq}^{\text{c}*}(s) = \boldsymbol{Y}_{\text{vir}}(s)\Delta\boldsymbol{u}_{\text{g}dq}(s) + \boldsymbol{Y}_{\text{vir}}(s)\boldsymbol{U}_{\text{g}dq0}^{\text{sync}}(s)\Delta\omega_{\text{sync}}(s) \tag{30}$$

where $\boldsymbol{Y}_{\text{vir}}(s)$ denotes the virtual admittance matrix. According to Fig. 3, the modulation ratio in the sync frame, i.e., $\Delta\boldsymbol{m}_{dq}(s)$, can be represented as:

$$\begin{aligned} \Delta\boldsymbol{m}_{dq}(s) &= -\boldsymbol{M}_{dq0}^{\text{sync}}(s)\Delta\omega_{\text{sync}}(s) - \\ &\frac{e^{-sT_{\text{s}}}}{V_{\text{dc0}}}\left[\boldsymbol{H}_{\text{cc}}(s)\left(\Delta\boldsymbol{i}_{\text{g}dq}^{\text{c}*}(s) - \Delta\boldsymbol{i}_{\text{g}dq}^{\text{c}}(s)\right) - \boldsymbol{K}_{\text{d}}\Delta\boldsymbol{i}_{\text{g}dq}^{\text{c}}(s)\right] \\ &= \boldsymbol{T}_{\text{cc1}}(s)\Delta\boldsymbol{i}_{\text{g}dq}^{\text{c}*}(s) + \boldsymbol{T}_{\text{cc2}}(s)\Delta\boldsymbol{i}_{\text{g}dq}(s) + \boldsymbol{T}_{\text{cc3}}(s)\Delta\omega_{\text{sync}}(s) \end{aligned} \tag{31}$$

where $\boldsymbol{H}_{\text{cc}}(s) = \text{diag}\left(H_{\text{cc}}, H_{\text{cc}}\right)$, $\boldsymbol{K}_{\text{d}} = \begin{bmatrix} 0 & -\omega_1 L_{\text{f}} \\ \omega_1 L_{\text{f}} & 0 \end{bmatrix}$. And $\boldsymbol{T}_{\text{cc1}}(s) = -\frac{e^{-sT_{\text{s}}}}{V_{\text{dc0}}}\boldsymbol{H}_{\text{cc}}$, $\boldsymbol{T}_{\text{cc2}}(s) = \frac{e^{-sT_{\text{s}}}}{V_{\text{dc0}}}\left(\boldsymbol{H}_{\text{cc}} + \boldsymbol{K}_{\text{d}}\right)$, $\boldsymbol{T}_{\text{cc3}}(s) = \frac{e^{-sT_{\text{s}}}}{V_{\text{dc0}}}\left(\boldsymbol{H}_{\text{cc}}\boldsymbol{I}_{\text{g}dq0}^{\text{sync}} + \boldsymbol{K}_{\text{d}}\boldsymbol{I}_{\text{g}dq0}^{\text{sync}}\right) - \boldsymbol{M}_{\text{g}dq0}^{\text{sync}}$. Then, converter-side voltage $\Delta\boldsymbol{u}_{\text{c}dq}(s)$ can be modeled as:

$$\Delta\boldsymbol{u}_{\text{c}dq}(s) = \Delta\boldsymbol{m}_{dq}(s)V_{\text{dc0}} + \boldsymbol{m}_{dq0}\Delta u_{\text{dc}}(s) \tag{32}$$

Based on (24), (25), (27)~(32), the converter's small-signal model under typical GFL/GFM control can be obtained as (33). By eliminating variables $\Delta\boldsymbol{u}_{\text{c}dq}$ and $\Delta\boldsymbol{i}_{\text{g}dq}^{\text{c}*}$, the relation between input $\left(\Delta\boldsymbol{u}_{\text{g}dq}, \Delta u_{\text{dc}}, \Delta f_{\text{sync}}\right)$ and output $\left(\Delta\boldsymbol{i}_{\text{g}dq}, \Delta i_{\text{dc}}, \Delta P_{\text{sync}}\right)$ can be established, i.e., four-port EIM shown in (9).

$$\begin{aligned}
\Delta\boldsymbol{u}_{cdq} &= G_1\left(\Delta\boldsymbol{i}_{gdq},\Delta\boldsymbol{u}_{gdq}\right)\\
\Delta i_{dc} &= G_2\left(\Delta\boldsymbol{i}_{gdq},\Delta\boldsymbol{u}_{cdq},\Delta u_{dc}\right)\\
\Delta P_{sync} &= \underbrace{G_3\left(\Delta\boldsymbol{u}_{gdq},\Delta f_{sync}\right)}_{GFL\ control}\ or\ \underbrace{G_4\left(\Delta\boldsymbol{i}_{gdq},\Delta\boldsymbol{u}_{gdq}\right)}_{GFM\ control}\\
\Delta\boldsymbol{i}^{c*}_{gdq} &= \underbrace{G_5\left(\Delta\boldsymbol{i}_{gdq},\Delta\boldsymbol{u}_{gdq},\Delta u_{dc}\right)}_{GFL\ control}\ or\ \underbrace{G_6\left(\Delta\boldsymbol{u}_{gdq},\Delta f_{sync}\right)}_{GFM\ control}\\
\Delta\boldsymbol{u}_{cdq} &= G_7\left(\Delta\boldsymbol{i}^{c*}_{gdq},\Delta\boldsymbol{i}_{gdq},\Delta f_{sync},\Delta u_{dc}\right)
\end{aligned} \tag{34}$$

## Appendix B
## Per Unit System Transformation

For the four-port model of GFL-type VSC, $\Delta P_{sync}$ is actually the $q$-axis voltage $\Delta u^{c}_{gq}$ and thus the rated voltage should be selected as the base value. Its base vector $\boldsymbol{I}_b$ and $\boldsymbol{U}_b$ is like:

$$\boldsymbol{I}_b = diag\left(\boldsymbol{I}^{ac}_b, I^{dc}_b, V^{ac}_b\right), \boldsymbol{U}_b = diag\left(\boldsymbol{V}^{ac}_b, V^{dc}_b, \omega_b\right) \tag{35}$$

For the four-port model of GFL-type VSC, $\Delta P_{sync}$ is exactly the active power and thus the rated power should be selected as the base value. Its base vector $\boldsymbol{I}_b$ and $\boldsymbol{U}_b$ is like:

$$\boldsymbol{I}_b = diag\left(\boldsymbol{I}^{ac}_b, I^{dc}_b, S_b\right), \boldsymbol{U}_b = diag\left(\boldsymbol{V}^{ac}_b, V^{dc}_b, \omega_b\right) \tag{36}$$

where $\boldsymbol{I}^{ac}_b$ and $I^{dc}_b$ are the ac-/dc- side rated current respectively; $\boldsymbol{V}^{ac}_b$ and $V^{dc}_b$ are the ac-/dc- side rated voltage respectively; $S_b$ denotes the rated power and $\omega_b$ is the rated frequency.

For $\boldsymbol{Z}_{net}(s)$ and $\boldsymbol{Y}_{con}(s)$ of the multi-VSCs system, the base vector $\boldsymbol{I}_b$ and $\boldsymbol{U}_b$ can be obtained by arranging (35) and (36) in a block diagonal form according to the port sequence.

## References

[1] "High Penetration of Power Electronic Interfaced Power Sources (HPoPEIPS)," ENTSO-E, Technical, Mar. 2017. [Online]. Available: https://consultations.entsoe.eu/system-development/entso-e-connection-codes-implementation-guidance-d-3.

[2] The Migrate Project [EB/OL]. https://www.h2020-migrate.eu/downloads.html. 2021-12-30.

[3] D. B. Rathnayake *et al.*, "Grid Forming Inverter Modeling, Control, and Applications," in *IEEE Access*, vol. 9, pp. 114781-114807, 2021.

[4] H. Zhang, W. Xiang, W. Lin and J. Wen, "Grid Forming Converters in Renewable Energy Sources Dominated Power Grid: Control Strategy, Stability, Application, and Challenges," in *Journal of Modern Power Systems and Clean Energy*, vol. 9, no. 6, pp. 1239-1256, November 2021.

[5] Y. Li, Y. Gu and T. C. Green, "Revisiting Grid-Forming and Grid-Following Inverters: A Duality Theory," in *IEEE Transactions on Power Systems*, vol. 37, no. 6, pp. 4541-4554, Nov. 2022.

[6] L. Huang *et al.*, "Grid-Synchronization Stability Analysis and Loop Shaping for PLL-Based Power Converters With Different Reactive Power Control," in *IEEE Transactions on Smart Grid*, vol. 11, no. 1, pp. 501-516, Jan. 2020.

[7] L. Huang *et al.*, "Impacts of grid structure on PLL-synchronization stability of converter integrated power systems," online. https：// arxiv.org / abs / 1903.05489.

[8] Y. Qi, H. Deng, J. Fang and Y. Tang, "Synchronization Stability Analysis of Grid-Forming Inverter: A Black Box Methodology," in *IEEE Transactions on Industrial Electronics*, vol. 69, no. 12, pp. 13069-13078, Dec. 2022.

[9] C. Yang, L. Huang, H. Xin and P. Ju, "Placing Grid-Forming Converters to Enhance Small Signal Stability of PLL-Integrated Power Systems," in *IEEE Transactions on Power Systems*, vol. 36, no. 4, pp. 3563-3573, July 2021.

[10] Wilsun Xu, Zhenyu Huang, Yu Cui and Haizhen Wang, "Harmonic resonance mode analysis," in *IEEE Transactions on Power Delivery*, vol. 20, no. 2, pp. 1182-1190, April 2005.

[11] Y. Zhan, X. Xie, H. Liu, H. Liu and Y. Li, "Frequency-Domain Modal Analysis of the Oscillatory Stability of Power Systems With High-Penetration Renewables," in *IEEE Transactions on Sustainable Energy*, vol. 10, no. 3, pp. 1534-1543, July 2019.

[12] D. Yang and Y. Sun, "SISO Impedance-Based Stability Analysis for System-Level Small-Signal Stability Assessment of Large-Scale Power Electronics-Dominated Power Systems," in *IEEE Transactions on Sustainable Energy*, vol. 13, no. 1, pp. 537-550, Jan. 2022.

[13] Y. Zhu, Y. Gu, Y. Li and T. C. Green, "Participation Analysis in Impedance Models: The Grey-Box Approach for Power System Stability," in *IEEE Transactions on Power Systems*, vol. 37, no. 1, pp. 343-353, Jan. 2022.

[14] H. Zong, C. Zhang, X. Cai and M. Molinas, "Oscillation Propagation Analysis of Hybrid AC/DC Grids With High Penetration Renewables," in *IEEE Transactions on Power Systems*, vol. 37, no. 6, pp. 4761-4772, Nov. 2022.

[15] C. Zhang, H. Zong, X. Cai and M. Molinas, "On the Relation of Nodal Admittance- and Loop Gain-Model Based Frequency-Domain Modal Methods for Converters-Dominated Systems," in *IEEE Transactions on Power Systems*, to be published.

[16] A. Rygg, M. Molinas, C. Zhang and X. Cai, "A Modified Sequence-Domain Impedance Definition and Its Equivalence to the dq-Domain Impedance Definition for the Stability Analysis of AC Power Electronic Systems," in *IEEE Journal of Emerging and Selected Topics in Power Electronics*, vol. 4, no. 4, pp. 1383-1396, Dec. 2016.

[17] C. Zhang, X. Cai, M. Molinas and A. Rygg, "On the Impedance Modeling and Equivalence of AC/DC-Side Stability Analysis of a Grid-Tied Type-IV Wind Turbine System," in *IEEE Transactions on Energy Conversion*, vol. 34, no. 2, pp. 1000-1009, June 2019.

[18] S. Shah and L. Parsa, "Impedance Modeling of Three-Phase Voltage Source Converters in DQ, Sequence, and Phasor Domains," in *IEEE Transactions on Energy Conversion*, vol. 32, no. 3, pp. 1139-1150, Sept. 2017.

[19] H. Liu, X. Xie, X. Gao, H. Liu and Y. Li, "Stability Analysis of SSR in Multiple Wind Farms Connected to Series-Compensated Systems Using Impedance Network Model," in *IEEE Transactions on Power Systems*, vol. 33, no. 3, pp. 3118-3128, May 2018.

[20] Y. Gu, Y. Li, Y. Zhu and T. C. Green, "Impedance-Based Whole-System Modeling for a Composite Grid via Embedding of Frame Dynamics," in *IEEE Transactions on Power Systems*, vol. 36, no. 1, pp. 336-345, Jan. 2021.

[21] L. Orellana, L. Sainz, E. Prieto-Araujo and O. Gomis-Bellmunt, "Stability Assessment for Multi-Infeed Grid-Connected VSCs Modeled in the Admittance Matrix Form," in *IEEE Transactions on Circuits and Systems I: Regular Papers*, vol. 68, no. 9, pp. 3758-3771, Sept. 2021.

[22] C. Zhang, M. Molinas, A. Rygg and X. Cai, "Impedance-Based Analysis of Interconnected Power Electronics Systems: Impedance Network Modeling and Comparative Studies of Stability Criteria," in *IEEE Journal of Emerging and Selected Topics in Power Electronics*, vol. 8, no. 3, pp. 2520-2533, Sept. 2020.

[23] J. Pedra, L. Sainz and L. Monjo, "Three-Port Small Signal Admittance-Based Model of VSCs for Studies of Multi-Terminal HVDC Hybrid AC/DC Transmission Grids," in *IEEE Transactions on Power Systems*, vol. 36, no. 1, pp. 732-743, Jan. 2021.

[24] H. Zhang, M. Mehrabankhomartash, M. Saeedifard, Y. Meng, X. Wang and X. Wang, "Stability Analysis of a Grid-Tied Interlinking Converter System With the Hybrid AC/DC Admittance Model and Determinant-Based GNC," in *IEEE Transactions on Power Delivery*, vol. 37, no. 2, pp. 798-812, April 2022.

[25] J. J.Grainger and J. William D. Stevenson, *Power system analysis*. New York, NY, USA: McGraw-Hill, 1994.